\documentclass{article}
\usepackage{amsmath, epsfig,amssymb, multicol}
\newtheorem{lemma}{Lemma}
\newtheorem{theorem}{Theorem}

\newtheorem{corollary}{Corollary}

\renewcommand{\L}{{\cal L}}

\renewcommand{\Pr}{{\bf Pr}}
\newcommand{\one}{{\bf 1}}
\newcommand{\vol}{{\rm vol}}
\newcommand{\diam}{{\rm diam}}
\newcommand{\vs}{{{\rm V}^{\underline{s}}}}
\newcommand{\vrm}{{{\rm V}^{\underline{r-1}}}}

\title{High-ordered Random Walks and Generalized Laplacians on Hypergraphs}
\author{Linyuan Lu
\thanks{University of South Carolina, Columbia, SC 29208,
({\tt lu@math.sc.edu}). This author was supported in part by NSF
grant DMS 1000475. }
  \and Xing Peng
\thanks{University of South Carolina, Columbia, SC 29208,
({\tt pengx@mailbox.sc.edu}).This author was supported in part by
NSF grant  DMS 1000475. }}
\begin{document}
\maketitle
\begin{abstract}
  Despite of the extreme success of the spectral graph theory, there
  are relatively few papers applying spectral analysis to hypergraphs.
  Chung first introduced Laplacians for regular hypergraphs and showed
  some useful applications. Other researchers treated
  hypergraphs as weighted graphs and then studied the Laplacians of the corresponding weighted graphs. In this paper, we aim to unify
  these very different versions of Laplacians for hypergraphs. We
  introduce a set of Laplacians for hypergraphs through studying
  high-ordered random walks on hypergraphs.  We prove the   eigenvalues of
  these Laplacians can effectively control the mixing rate of high-ordered
  random walks, the generalized distances/diameters, and  the edge expansions.
\end{abstract}

\section{Introduction}
Many complex networks have richer structures than graphs can have.
Inherently they have  hypergraph structures: interconnections often
cross multiple nodes. Treating these networks as graphs causes a
loss of some structures. Nonetheless, it is still popular to use
graph tools to study these networks; one of them is the
Laplacian spectrum. Let $G$ be a graph on $n$ vertices.
The {\it Laplacian} $\L$ of  $G$ is the $(n\times n)$-matrix
$I-T^{-1/2}AT^{-1/2}$, where $A$ is the adjacency matrix and $T$ is the diagonal
matrix of degrees.  Let $\lambda_0, \lambda_1, \ldots, \lambda_{n-1}$ be the eigenvalues of $\L$, indexed in
non-decreasing order. It is known that $0\leq \lambda_i \leq 2$ for $0\leq i\leq n-1$.
If $G$ is connected, then $\lambda_1>0$. The
first nonzero Laplacian eigenvalue $\lambda_1$  is related to many
graph parameters, such as the mixing rate of random walks, the graph
diameter, the neighborhood expansion, the Cheeger constant, the
isoperimetric inequalities, expander graphs, quasi-random graphs,
etc \cite{alon1,alon2,fan1,fan3,fan4}.

In this paper, we define a set of Laplacians for hypergraphs.
Laplacians for regular hypergraphs was first introduced by Chung
\cite{fan2} in 1993 using homology approach. The first nonzero
Laplacian eigenvalue can be used to derive several useful
isoperimetric inequalities. It seems hard to extend Chung's
definition to general hypergraphs. Other researchers  treated a
hypergraph as a multi-edge graph and then defined its Laplacian to
be the Laplacian of the corresponding multi-edge graph. For example,
Rodr\'iguez \cite{rod} showed that the approach above had some
applications to bisections, the average minimal cut, the
isoperimetric number, the max-cut, the independence number, the
diameter etc.

What are ``right'' Laplacians for hypergraphs? To answer this
question, let us recall how the Laplacian was introduced in the
graph theory. One of the approaches is using geometric/homological analogue,
where the
Laplacian is defined as a self-joint operator on the functions over
vertices. Another approach is using random walks, where the
Laplacian is the symmetrization of the transition matrix of the random walk
on a graph. Chung \cite{fan1} took the first approach and defined her Laplacians
for regular hypergraphs. In this paper, we
take the second approach and define the Laplacians through
high-ordered random walks on hypergraphs.

 A high-ordered walk on a hypergraph $H$ can be roughly viewed as a sequence of
overlapped oriented edges $F_1, F_2,\ldots, F_k$. For $1\leq s\leq r-1$, we
say $F_1, F_2,\ldots, F_k$ is an $s$-walk if  $|F_i\cap F_{i+1}|=s$
for each $i$ in $\{1,2,3,\ldots,k-1\}$. The choice of $s$ enables us to define
a set of Laplacian matrices $\L^{(s)}$ for $H$. For $s=1$, our definition of
Laplacian $\L^{(1)}$ is the same as the definition in  \cite{rod}.
For $s=r-1$, while we restrict to regular hypergraphs, our
definition of Laplacian $\L^{(r-1)}$ is similar to Chung's
definition \cite{fan2}. We will discuss their relations in the last section.

In this paper, we show several applications of the Laplacians of hypergraphs,
such as the mixing rate of high-ordered random walks, the
generalized diameters, and the edge expansions. Our approach allows
users to select a ``right'' Laplacian to fit their special
need.

The rest of the paper is organized as  follows. In section 2, we
review and prove some useful results on the Laplacians of weighted
graphs and Eulerian directed graphs. The definition of  Laplacians
for hypergraphs will be  given in section 3. We will
prove some properties of the  Laplacians of hypergraphs in section 4,
and consider several applications in section 5.
In last section, we will comment on future
directions.

\section{Preliminary results}
In this section, we review some results on Laplacians of  weighted graphs and Eulerian directed graphs.
Those results will be
applied to the Laplacians of hypergraphs later on.

In this paper, we frequently switch domains from hypergraphs to weighted (undirected)
graphs, and/or to directed graphs. To reduce confusion, we use the following
conventions through this paper. We denote a weighted graph by $G$, a
directed graph by $D$, and a hypergraph by $H$. The set of vertices
is denoted by $V(G)$, $V(D)$, and $V(H)$, respectively. (Whenever it
is clear under the  context, we will write it as $V$ for short.) The
edge set is denoted by $E(G)$, $E(D)$, and $E(H)$, respectively. The
degrees $d_\ast$ and volumes $\vol(\ast)$ are defined separately for
the weighted graph $G$, for the directed graph $D$, and for the
hypergraph $H$. Readers are warned to interpret them   carefully under
the context.

For a positive integer $s$ and a vertex set $V$, let $\vs$ be the set of all (ordered)
$s$-tuples consisting of $s$ distinct elements in $V$.  Let
${V\choose
  s}$ be the set of all unordered (distinct) $s$-subset of $V$.

Let $\one$ be the row (or the column) vector with all entries of value $1$
 and $I$ be the identity matrix. For a row (or column) vector $f$, the norm
$\|f\|$ is always the $L_2$-norm of $f$.

\subsection{Laplacians of weighted graphs}
A {\it weighted graph} $G$  on the vertex set $V$ is an undirected
graph associated with  a weight function $w \colon V\times V \to
\mathbb R^{\geq 0}$ satisfying  $w(u, v) = w(v, u)$ for all $u$ and
$v$ in $V(G)$. Here we always assume $w(v,v)=0$ for every $v\in V$.

A simple graph can be viewed as a special weighted graph with weight
$1$ on all edges and $0$ otherwise. Many concepts of simple graphs
are naturally generalized to weighted graphs. If $w(u,v)>0$, then
$u$ and $v$ are {\it adjacent}, written as  $x\sim y$. The graph
distance $d(u,v)$ between two vertices $u$ and $v$ in $G$ is the
minimum integer $k$ such that there is a path
$u=v_0,v_1,\ldots,v_k=v$ in which $w(v_{i-1}, v_{i})>0$ for $1\leq
i\leq k$. If no such $k$ exists, then we let $d(u,v)=\infty$. If the
distance $d(u,v)$ is finite for every pair $(u,v)$,  then $G$ is
{\it connected}. For a connected weighted graph $G$, the {\it
diameter} (denoted by $\textrm{diam}(G)$) is the smallest value of
$d(u,v)$ among all pairs of vertices $(u,v)$.

The {\it adjacency matrix} $A$ of $G$ is defined as the matrix of
weights, i.e., $A(x,y)=w(x,y)$  for all $x$ and $y$ in $V$. The {\it
degree} $d_x$ of a vertex $x$ is $\sum_y w(x,y)$.  Let $T$ be the
diagonal matrix of degrees in $G$. The {\it Laplacian} $\L$ is the
matrix $I-T^{-1/2}AT^{-1/2}$. Let $\lambda_0, \lambda_1, \ldots,
\lambda_{n-1}$ be the eigenvalues of $\L$, indexed in the
non-decreasing order. It is known \cite{fan4}  that $0\leq \lambda_i
\leq 2$ for $0\leq i\leq n-1$. If $G$ is connected, then
$\lambda_1>0$.

From now on, we assume $G$ is connected. The
first non-trivial  Laplacian eigenvalue $\lambda_1$ is the
most useful one. It can be written in terms of  the Rayleigh quotient as follows
(see \cite{fan4})
\begin{equation}
  \label{eq:2}
  \lambda_1=\inf_{f\perp T\one} \frac{\sum_{x \sim y} (f(x)-f(y))^2w(x,y)}{\sum_x f(x)^2 d_x}.
\end{equation}
Here the infimum is taken over all functions $f\colon V\to \mathbb
R$ which is orthogonal to  the degree vector $\one
T=(d_1,d_2,\ldots, d_n)$.
Similarly, the largest Laplacian eigenvalue
$\lambda_{n-1}$ can be defined in terms of the Rayleigh quotient
as follows
\begin{equation}
  \label{eq:3}
  \lambda_{n-1}=\sup_{ f\perp T\one} \frac{\sum_{x \sim y} (f(x)-f(y))^2w(x,y)}{\sum_x f(x)^2 d_x}.
\end{equation}

Note that scaling the weights by a constant factor will not affect
the Laplacian. A weighted graph $G$ is {\it complete} if $w(u,v)=c$
for some constant $c$ such that $c>0$, independent of the choice of
$(u,v)$ with $u\not=v$. We say $G$ is {\it bipartite} if there is a
partition $V=L\cup R$ such that $w(x,y)=0$ for all $x,y\in L$ and
all $x,y\in R$.

We have the following facts (see \cite{fan4}).
\begin{enumerate}
\item  $0 \leq \lambda_i \leq 2$ for each $0 \leq i \leq n-1$.
\item The number of $0$ eigenvalues equals the number of connected components
in $G$. If $G$ is connected, then $\lambda_1>0$.
\item $\lambda_{n-1}=2$ if and only if  $G$ has a connected component
which is a bipartite weighted subgraph.
\item $\lambda_{n-1}=\lambda_1$ if and only if $G$ is a complete weighted graph.
\end{enumerate}

It turns out that $\lambda_1$ and $\lambda_{n-1}$ are related to
many graph parameters, such as the mixing rate of random walks, the
diameter, the edge expansions, and the isoperimetric inequalities.

A random walk on a weighted graph $G$ is a sequence of  vertices
$v_0, v_1,\ldots, v_k$ such that the conditional probability
$\Pr(v_{i+1}=v\mid v_i=u)=w(u,v)/d_u$ for $0\leq i\leq k-1$. A {\it
vertex probability distribution} is a map $f\colon V\to \mathbb R$
such that  $f(v)\geq 0$ for each $v$ in $G$ and $\sum_{v\in
V}f(v)=1$. It is convenient to write a vertex probability
distribution into a row vector.  A random walk maps a vertex
probability distribution to a vertex probability distribution
through multiplying from right a transition matrix $P$, where
$P(u,v)=w(u,v)/d_u$ for each pair of vertices $u$ and $v$. We can
write $P=T^{-1}A=T^{-1/2}(I-\L)T^{1/2}$. The second largest
eigenvalue $\bar \lambda(P)$, denoted by $\bar \lambda$ for short,
 is $\max\{|1-\lambda_1|, |1-\lambda_{n-1}|\}$. Let
$\pi(u)=d_u/\vol(G)$ for each vertex $u$ in $G$. Observe $\pi$
 is the stationary distribution of the random
walk, i.e.,  $\pi P=\pi$. A random walk is {\it mixing} if ${\lim_{i
\rightarrow \infty}} f_0 P^i=\pi$ for any initial vertex probability distribution $f_0$.
 It is known that a random walk is
always mixing if $G$ is connected and not a bipartite graph. To
overcome the difficulty resulted from being a bipartite graph (where
$\lambda_{n-1}=2$),
for $0\leq \alpha \leq 1$, we consider an $\alpha$-lazy random walk, whose
transition matrix $P_\alpha$ is given by $P_\alpha(u,u)=\alpha$ for each $u$
and $P_\alpha(u,v)=(1-\alpha)
w(u,v)/d_u$ for each pair of vertices $u$ and $v$ with $u\not=v$.
  Note that the transition matrix is
$$P_\alpha=\alpha I + (1-\alpha) T^{-1}A=T^{-1/2}(
I-(1-\alpha)\L)T^{1/2}.$$
Let $L_{\alpha}=T^{1/2}P_{\alpha}T^{-1/2}= I-(1-\alpha)\L$
and $\bar\lambda_\alpha=\max\{|1-(1-\alpha)\lambda_1|, |1-(1-\alpha)\lambda_{n-1}|\}$.
Since $L_{\alpha}$ is a symmetric matrix, we have
$$
\bar \lambda_{\alpha}=\underset{u \perp  T^{1/2} \one}\max
\frac{\|L_{\alpha}u \|}{\|u\|}.
$$

 It turns out that the mixing rate of an $\alpha$-lazy random walk
is determined by $\bar \lambda_{\alpha}$.
\begin{theorem}\label{thm1}
For $0 \leq \alpha \leq 1$, the vertex probability distribution
$f_k$ of the $\alpha$-lazy random walk at time $k$ converges to the
stationary distribution $\pi$ in probability. In particular, we have
$$\| (f_k-\pi)T^{-1/2} \| \leq \bar \lambda^k \|(f_0-\pi)T^{-1/2}\|.$$
Here $f_0$ is the initial vertex probability distribution.
\end{theorem}
{\bf Proof:} Notice that $f_k=f_0P_{\alpha}^k$ and
$(f_0-\pi)T^{-1/2} \perp \one T^{1/2} $. We have
\begin{eqnarray*}\hspace*{1.3in}
\| (f_k-\pi)T^{-1/2} \|& = & \|(f_0P_{\alpha}^k -\pi P_{\alpha}^k)T^{-1/2}\|\\
             &=& \| (f_0-\pi)P_{\alpha}^kT^{-1/2}\|\\
             &=& \| (f_0-\pi) T^{-1/2}L_{\alpha}^k\| \\
             & \leq & \bar \lambda_{\alpha}^k
             \|(f_0-\pi)T^{-1/2}\|. \hspace*{1.3in}\hfill \square
           \end{eqnarray*}

For each subset $X$ of  $V(G)$, the {\it volume} $\vol(X)$ is $\sum_{x\in
X}d_x$. If $X=V(G)$, then we write $\vol(G)$ instead of $\vol(V(G))$. We
have
$$\vol(G)=\sum_{i=1}^n d_i=2\sum_{u\sim v} w(u,v).$$
If $\bar X$ is the complement set of $X$,  then have $\vol(\bar X)=\vol(G)-\vol(X)$.
For any two subsets $X$ and $Y$ of $V(G)$, the {\it distance} $d(X,Y)$ between $X$ and $Y$
is $\min\{d(x,y): x\in X, y \in Y\}$.

\begin{theorem}[See \cite{fan1,fan4}]\label{thm2}
In a weighted graph $G$, for $X,Y \subseteq V(G)$ with distance at least 2, we have
$$
d(X,Y) \leq \left \lceil \frac{\log \sqrt {\frac{\vol(\bar
X)\vol(\bar Y)}{\vol(X)\vol(Y)}}}
  {\log
  \frac{\lambda_{n-1}+\lambda_1}{\lambda_{n-1}-\lambda_1}}\right\rceil.
$$
\end{theorem}

A special case of Theorem \ref{thm2} is that  both $X$ and $Y$ are
single vertices,  which gives an upper bound on the diameter of $G$.
\begin{corollary}[See \cite{fan4}]\label{cor1}
If $G$ is not a complete weighted graph, then we have
$$
\textrm{diam}(G) \leq   \left \lceil
\frac{\log (\vol(G)/\delta)}
  {\log
  \frac{\lambda_{n-1}+\lambda_1}{\lambda_{n-1}-\lambda_1}}\right\rceil,
$$
where $\delta$  is the minimum degree of $G$.
\end{corollary}
For $X,Y \subseteq V(G)$, let $E(X,Y)$ be the set of edges between $X$ and $Y$.
Namely, we have
$$
E(X,Y)=\{(u,v): u \in X, v \in Y \ \textrm{and}\ uv\in E(G) \}.
$$

We have the following theorem.
\begin{theorem}[See \cite{fan1,fan4}]\label{thm3}
If $X$ and $Y$ are two subsets of $V(G)$, then we have
$$
\left||E(X,Y)|-\frac{\vol(X) \vol(Y)}{\vol(G)}\right| \leq \bar
\lambda \frac{\sqrt{\vol(X)\vol(Y)\vol(\bar X) \vol(\bar
Y)}}{\vol(G)}.
$$
\end{theorem}

\subsection{Laplacians of Eulerian directed graphs}
The Laplacian of a general directed graph was introduced by Chung
\cite{fan5, fan6}. The theory is considerably more complicated than
the one for undirected graphs, but when we  consider a special class
of directed graphs --- Eulerian directed graphs,  it turns out
to be quite neat.

Let $D$ be a directed graph with the vertex set $V(D)$ and the edge set $E(D)$.
A directed edge from $x$ to $y$ is denoted by an ordered pair $(x,y)$ or $x\to y$.
The {\it out-neighborhood} $\Gamma^+(x)$ of a vertex $x$ in $D$ is
the set $\{y\colon (x,y)\in E(D)\}$. The {\it out-degree} $d^+_x$  is $|\Gamma^+(x)|$.
Similarly, the {\it in-neighborhood} $\Gamma^-(x)$ is  $\{y\colon (y,x)\in E(D)\}$,
and  the {\it in-degree} $d^-_x$  is $|\Gamma^-(x)|$.
A directed graph $D$ is {\it Eulerian} if
$d^+_x=d^-_x$ for every vertex $x$.  In this case, we simply write $d_x=d^+_x=d^-_x$
for each $x$. For
a vertex subset $S$, the {\it volume} of $S$, denoted by $\vol(S)$, is $\sum_{x\in S}d_x$.
In particular, we write $\vol(D)=\sum_{x\in V}d_x$.

Eulerian directed graphs have many good properties. For example, a
Eulerian directed graph is strongly connected if and only if it is
weakly connected.

The {\it adjacency matrix} of $D$ is a square matrix $A$ satisfying
 $A(x,y)=1$ if $(x,y) \in E(D)$ and 0 otherwise.  Let $T$ be the
diagonal  matrix with $T(x,x)=d_x$ for each $x\in V(D)$.
Let $\vec
\L=I-T^{-1/2}AT^{-1/2}$, i.e.,
\begin{equation}
  \label{eq:ls}
 \vec\L(x,y)=\left\{
   \begin{array}{ll}
     1 & \mbox{ if } x=y;\\
     -\frac{1}{\sqrt{d_xd_y}} & \mbox{ if } x\to y;\\
0 & \mbox{ otherwise.}
   \end{array}
\right.
\end{equation}

Note that $\vec \L$ is not symmetric. We define the Laplacian $\L$
of $D$ to be the symmetrization  of $\vec \L$, that is
$$\L=\frac{\vec \L + \vec \L'}{2}.$$
Since $\L$ is symmetric, its eigenvalues are real and can be listed
as $\lambda_0,\lambda_1,\ldots,\lambda_{n-1}$ in     the
non-decreasing order. Note that $\lambda_1$ can also  be written in
terms of Raleigh quotient (see \cite{fan5}) as follows
\begin{equation}
  \label{eq:rl2}
  \lambda_1=\inf_{ f\perp T\one} \frac{\sum_{x \to y} (f(x)-f(y))^2}{2\sum_x f(x)^2 d_x}.
\end{equation}

Chung \cite{fan6} proved a general theorem on the relationship
between $\lambda_1$ and the
diameter. After restricting to Eulerian
directed graphs, it can be stated as follows.

\begin{theorem} [See \cite{fan6}]\label{thm4}
Suppose $D$ is  a connected Eulerian directed graph, then the
diameter of $D$ (denoted by $diam(D)$) satisfies
$$ diam(D) \leq \left \lfloor \frac{2\log(\vol(G)/\delta)}{\log
\frac{2}{2-\lambda_1}} \right \rfloor +1,
$$
where $\lambda_1$ is the first non-trivial eigenvalue of the
Laplacian, and $\delta$ is the minimum degree $\min\{d_x\mid x\in V(D)\}$.
\end{theorem}

The main idea in the proof of the theorem above is using
$\alpha$-lazy random walks on $D$. A random walk on a Eulerian
directed graph $D$ is
 a sequence of  vertices $v_0, v_1,\ldots, v_k$
such that for $0\leq i\leq k-1$, the conditional probability
$\Pr(v_{i+1}=v\mid v_i=u)$ equals $1/d_u$ for each $v\in \Gamma^+(u)$
and $0$ otherwise. For $0\leq \alpha\leq 1$, the
$\alpha$-lazy random walk is defined similarly. The transition
matrix  $P_\alpha$ of the $\alpha$-lazy random walk satisfies
$$P_\alpha=\alpha I + (1-\alpha) T^{-1}A=T^{-1/2}(I-(1-\alpha)\vec\L)T^{1/2}.$$

Chung \cite{fan5} considered only $1/2$-lazy random walks. Here we
prove some results on $\alpha$-lazy random walks for $\alpha\in
[0,1)$.

Let $\pi(u)=d_u/ \vol(D)$ for each $u \in V(D)$. Note that $\pi$ is the stationary distribution, i.e.
$\pi P_\alpha=\pi$. Let $L_\alpha=\alpha I + (1-\alpha)
T^{-1/2}AT^{-1/2}=I-(1-\alpha)\vec\L=T^{1/2}P_\alpha T^{-1/2}$. The
key observation is that there is a unit-vector $\phi_0$ such that
$\phi_0$ is both a row eigenvector and a column eigenvector of $L_{\alpha}$
for the largest eigenvalue $1$. Here let $\phi_0= \one
T^{1/2}/{\vol(D)} = \frac{1}{\vol(G)}(\sqrt{d_1},\ldots,
\sqrt{d_n})$.  We have
$$\phi_0 L_\alpha = \phi_0  \mbox{ and } L_\alpha\phi_0'=\phi_0'.$$
Let $\phi_0^{\perp}$ be the orthogonal complement of $\phi_0$ in
$R^n$. It is easy to check $L_\alpha$ maps $\phi_0^{\perp}$ to
$\phi_0^{\perp}$. Let $\sigma_\alpha$ be the spectral norm of
$L_\alpha$ when restricting to $\phi_0^{\perp}$. An equivalent
definition of $\sigma_\alpha$ is the second largest singular
value of $L_\alpha$, i.e.,
$$
\sigma_{\alpha}=\underset{f \perp \phi_0'}\max
\frac{\|L_{\alpha}f\| }{\|f\|}.
$$
\begin{lemma}\label{lemma1}
We have the following properties for $\sigma_\alpha$.
\begin{enumerate}
\item For every $\beta \in \phi_0^{\perp}$, we have $ \|L_\alpha\beta\| \leq \sigma_\alpha \|\beta\|$.

\item $(1-\lambda_1)^2\leq \sigma_0^2 \leq 1.$
\item $\sigma^2_\alpha\leq \alpha^2 + 2\alpha (1-\alpha)\lambda_1 + (1-\alpha)^2 \sigma^2_0.$
\end{enumerate}
\end{lemma}
{\bf Proof:} Item 1 is from the definition of $\sigma_\alpha$.
Since the largest eigenvalue of
$L_{\alpha}$ is 1, we have $\sigma_{\alpha} \leq 1$. In particular,
$\sigma_0^2 \leq 1$. Note that $L_0=T^{-1/2}AT^{-1/2}$. Let
$f=gT^{1/2} $. It follows that
 $$\sigma_0^2=\underset{f \perp \phi_0'}\sup
\frac{\|L_{0}f\|^2}{\|f\|^2}=\underset{g \perp  T \one }\sup
\frac{g'A'T^{-1}Ag }{g'Tg}.$$

Choose $g\in (T\one)^\perp$ such that the Raleigh quotient (\ref{eq:rl2}) reaches its minimum at $g$, i.e.,
$$\lambda_1=\frac{\sum_{x \to y} (g(x)-g(y))^2}{2\sum_x g(x)^2 d_x}.$$

 We have
\begin{eqnarray*}
\frac{g'A'T^{-1}Ag }{g'T g}&=&\frac{\sum_x\frac{1}{d_x}\left (
\sum_{y\in \Gamma^+(x)} g(y)\right)^2 }{\sum_x d_xg(x)^2}\\
&=&\frac{\sum_x d_xg(x)^2 \sum_x\frac{1}{d_x}\left (
\sum_{y\in \Gamma^+(x)} g(y)\right)^2}{(\sum_x d_xg(x)^2)^2}\\
& \geq & \frac{\left(\sum_x g(x) \sum_{y\in \Gamma^+(x)}g(y)\right)^2
}{(\sum_x d_x g(x)^2)^2}\\
&=& \left( \frac{\sum_x g(x) \sum_{y\in \Gamma^+(x)}g(y)}{\sum_x d_x
g(x)^2}    \right)^2\\
&=&(1-\lambda_1)^2.
\end{eqnarray*}
In the last step, we use the following argument.
\begin{eqnarray*}
\frac{\sum_x g(x) \sum_{y\in \Gamma^+(x)}g(y)}{\sum_x d_x g(x)^2}
&=&\frac{\frac{1}{2}\sum_{x \rightarrow
y}\left(g(x)^2+g(y)^2-(g(x)-g(y))^2\right) }{\sum_x d_x g(x)^2}\\
&=&1-\frac{\sum_{x \rightarrow y}\left( g(x)-g(y) \right)^2}{2\sum_x
d_xg(x)^2}\\
&=& 1-\lambda_1.
\end{eqnarray*}
Since $\sigma_0$ is the maximum over all $g \perp T \one$, we get
$(1-\lambda_1)^2 \leq \sigma_0^2$.

For item 3, we have
\begin{eqnarray*}\hspace*{0.5in}
  \sigma_\alpha^2 &=& \sup_{f \perp \phi_0'}
\frac{\|L_{\alpha}f\|^2}{\|f\|^2}\\
&=&\underset{g \perp  T \one }\sup
\frac{g'P_\alpha'TP_\alpha g }{g'Tg}\\
&\leq&\alpha^2 + \alpha(1-\alpha)\sup_{g \perp  T \one}\frac{g'(A+A')g }{g'Tg}
+(1-\alpha)^2 \sup_{g \perp  T \one}\frac{g'A'T^{-1}Ag }{g'Tg}\\
&=& \alpha^2 + 2\alpha(1-\alpha)(1-\lambda_1) + (1-\alpha)^2\sigma_0^2.
\hspace*{1.7in} \square
\end{eqnarray*}

\begin{theorem}\label{thm5}
  For $0<\alpha<1$, the vertex probability distribution $f_k$ of the
  $\alpha$-lazy random walk on a Eulerian directed graph $D$ at time
  $k$ converges to the stationary distribution $\pi$ in probability. In
  particular, we have
$$\| (f_k-\pi)T^{-1/2} \| \leq \sigma_\alpha^k  \|(f_0-\pi)T^{-1/2}\|.$$
Here $f_0$ is the initial vertex probability distribution.
\end{theorem}
The proof is omitted since it is very similar to the proof of
Theorem \ref{thm1}. Notice that when $0 < \alpha < 1$, we have
$\sigma_{\alpha} < 1$ by Lemma \ref{lemma1}. We have  the
$\alpha$-lazy random converges to the stationary distribution
exponentially fast.


 For two vertex subsets $X$ and $Y$ of $V(D)$, let $E(X,Y)$ be the number of directed
edges from $X$ to $Y$, i.e.,
 $E(X,Y)=\{(u,v): u \in X \ \mbox{and}  \  v \in Y\}$.
We have the following theorem on the edge expansions in Eulerian
directed graphs.

\begin{theorem}\label{thm6}
If $X$ and  $Y$ are two subsets of the vertex set $V$ of a Eulerian
directed graph $D$, then we have
$$
\left||E(X,Y)|-\frac{\vol(X) \vol(Y)}{\vol(D)}\right| \leq \sigma_0
\frac{\sqrt{\vol(X)\vol(Y)\vol(\bar X) \vol(\bar Y)}}{\vol(D)}.
$$
\end{theorem}
{\bf Proof:} Let $\one_X$  be the indicator variable of $X$, i.e.,
$\one_X(u)=1$ if $u \in X$ and 0 otherwise.  We define $\one_Y$
similarly. Assume $\one_X T^{1/2}=a_0 \phi_0+a_1 \phi_1$ and
$\one_YT^{1/2}=b_0 \phi_0 + b_1 \phi_2 $, where $\phi_1,\phi_2 \in
\phi_0^{\perp}$ and are unit vectors. Since $\phi_0$ is a unit
vector,
 we have
\begin{equation}
  \label{eq:a0}
a_0= \langle \one_XT^{1/2}, \phi_0 \rangle=\frac{\vol(X)}{\sqrt{\vol(D)}}
\end{equation}
 and
 \begin{equation}
   \label{eq:a0b0}
a_0^2+a_1^2=\langle \one_X T^{1/2}, \one_X T^{1/2} \rangle=\vol(X).
 \end{equation}
Thus
\begin{equation}
  \label{eq:b0}
 a_1=\sqrt{\vol(X)\vol(\bar X)/\vol(D)}.
\end{equation}

Similarly, we get
\begin{eqnarray}
  \label{eq:a1}
b_0&=&\frac{\vol(Y)}{\sqrt{\vol(D)}}; \\
\label{eq:b1}
b_1&=&\sqrt{\vol(Y)\vol(\bar Y)/\vol(D)}.
\end{eqnarray}
It follows that
\begin{eqnarray*}
\left||E(X,Y)|-\frac{\vol(X) \vol(Y)}{\vol(D)}\right|
&=& \left |
\one_X T^{1/2} (L_0- \phi_0' \phi_0 ) (\one_Y T^{1/2})'
 \right |\\
 &=&\left |(a_0\phi_0+a_1\phi_1) (L_0- \phi_0' \phi_0 ) (b_0\phi_0+b_1\phi_2)'
 \right |\\
 &=& \left| a_1b_1 \phi_1 L_0 \phi_2'   \right |\\
 & \leq & |a_1b_1| \| \phi_1\| \| L_0 \phi_2' \| \\
 & \leq & |a_1b_1| \sigma_0 \\
 &=&\sigma_0
\frac{\sqrt{\vol(X)\vol(Y)\vol(\bar X) \vol(\bar Y)}}{\vol(D)}.
\end{eqnarray*}
The proof of this theorem is completed. \hfill $\square$

If we use $\bar\lambda$ instead of $\sigma_0$, then we get a weaker
theorem on the edge expansions. The proof will be omitted since it is very similar to
the proof of Theorem \ref{thm6}.

\begin{theorem}\label{thm7}
  Let $D$ be a Eulerian directed graph.  If $X$ and $Y$ are two
  subsets of $V(D)$, then we have
$$
\left|\frac{|E(X,Y)|+ |E(Y,X)|}{2}-\frac{\vol(X) \vol(Y)}{\vol(D)}\right| \leq \bar \lambda
\frac{\sqrt{\vol(X)\vol(Y)\vol(\bar X) \vol(\bar Y)}}{\vol(D)}.
$$
\end{theorem}

For $X,Y \subseteq V(D)$, let
$d(X,Y)=\min \{d(u,v): u \in X \ \textrm{and} \ v \in Y\}$. We have the following upper bound on
$d(X,Y)$.

\begin{theorem}\label{thm8}
Suppose $D$ is a connected Eulerian directed graph. For $X,Y \subseteq
V(D)$ and $0 \leq \alpha <1$, we have
$$
d(X,Y) \leq \left \lfloor \frac{\log \sqrt {\frac{\vol(\bar
X)\vol(\bar Y)}{\vol(X)\vol(Y)}}}
  {\log \sigma_\alpha}\right\rfloor+1.
$$
In particular, for $0\leq \alpha <1$, the diameter of $D$ satisfies
$$
\diam(D) \leq \left \lceil \frac{\log (\vol(D)/\delta)}
  {\log \sigma_\alpha}\right\rceil,
$$
 where $\delta=\min\{d_x\colon x\in V\}$.
\end{theorem}
{\bf Remark:}
From lemma \ref{lemma1}, we have
$$\sigma_\alpha^2 \leq \alpha^2+2\alpha(1-\alpha)\lambda_1
+(1-\alpha)^2 \sigma_0^2.$$
We can choose $\alpha$ to minimize $\sigma_\alpha$.
If $\lambda_1\leq 1-\sigma_0^2$, then we choose $\alpha=0$
and get $\sigma_\alpha=\sigma_0$; if $\lambda_1> 1-\sigma_0^2$,
then we choose $\alpha=\frac{\lambda_1+\sigma_0^2-1}{2\lambda_1+\sigma_0^2-1}$
and get $\sigma_\alpha^2\leq 1-\frac{\lambda_1^2}{2\lambda_1+\sigma_0^2-1}$.
Combining two cases, we have
\begin{equation}
  \label{eq:minalpha}
\min_{0\leq \alpha <1}\{\sigma_\alpha\} \leq \left\{
  \begin{array}{ll}
    \sigma_0 & \mbox{ if } \lambda_1\leq 1-\sigma_0^2; \\
   \sqrt{1-\frac{\lambda_1^2}{2\lambda_1+\sigma_0^2-1}}
&  \mbox{ otherwise.}
  \end{array}
\right.
\end{equation}
It is easy to check
$$\min_{0\leq \alpha <1}\{\sigma_\alpha\}\leq \sqrt{1-\frac{\lambda_1}{2}}.$$
Here the inequality is strict if $\sigma_0<1$. We have
$$
d(X,Y) \leq \left \lfloor \frac{\log \frac{\vol(\bar
X)\vol(\bar Y)}{\vol(X)\vol(Y)}}
  {\log \frac{2}{2-\lambda_1}}\right\rfloor +1.$$
Theorem \ref{thm8}
is stronger  than Theorem \ref{thm4} in general.

{\bf Proof:} Similar to the proof of Theorem \ref{thm6},
let $\one_X $ and $\one_Y$ be the indicator
functions of $X$ and $Y$, respectively. We have
\begin{eqnarray*}
 \one_X T^{1/2}&=&a_0 \phi_0+a_1 \phi_1,\\
 \one_YT^{1/2}&=&b_0 \phi_0 + b_1 \phi_2,
\end{eqnarray*}
where $\phi_1,\phi_2
\in \phi_0^{\perp}$ and are unit vectors and $a_0,b_0, a_1, b_1$
are given by equations (\ref{eq:a0})-(\ref{eq:b1}).

Let $k=\left\lfloor \frac{\log \sqrt {\frac{\vol(\bar
X)\vol(\bar Y)}{\vol(X)\vol(Y)}}}
  {\log \sigma_\alpha}\right\rfloor +1$. We have
$$(\one_X
T^{1/2})L_{\alpha}^k (\one_Y T^{1/2})'\geq a_0b_0+\sigma_\alpha^k a_1b_1> 0.$$

Thus there is a directed path starting from some vertex in $X$ and ending at some vertex in $Y$,
that is $d(X,Y) \leq k$.

For the diameter result, we choose $X=\{x\}$ and $Y=\{y\}$.
Note that $\vol(X)=d_x\geq \delta$, $\vol(Y)=d_y\geq \delta$,
$\vol(\bar X)<\vol(G)$, and $\vol(\bar y)<\vol(G)$.
The result follows.
\hfill $\square$

%
%
%
%
%

\section{Definition of the $s$-th Laplacian}
\label{section3}
\begin{figure}[hbt]
 \centerline{\psfig{figure=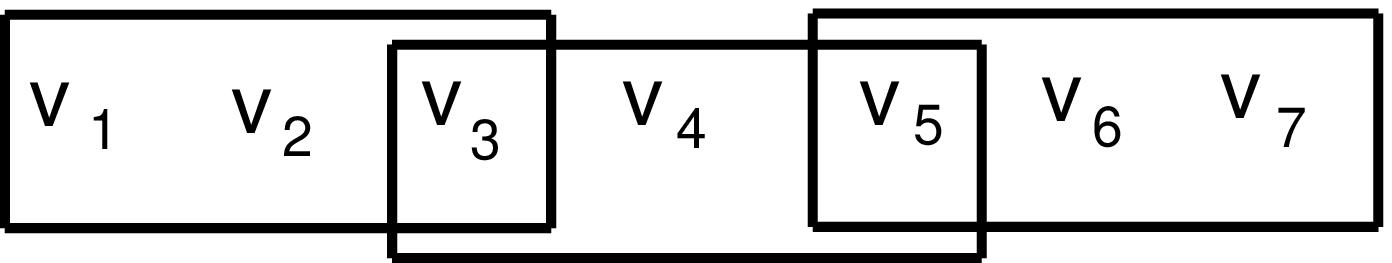, width=0.32\textwidth} \hfil \psfig{figure=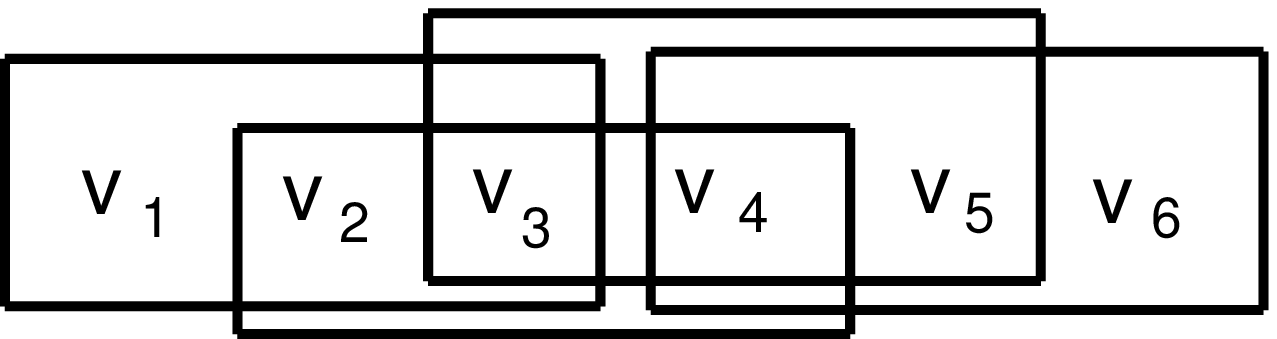,
 width=0.32\textwidth}  \hfil \psfig{figure=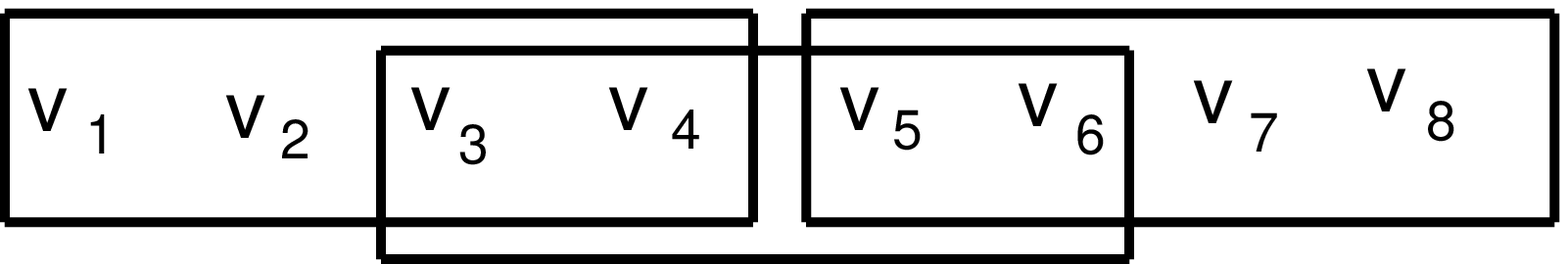,width=0.32\textwidth}}
 \centerline{A 1-walk in a 3-graph \hspace*{1cm} \hfil A 2-walk in a 3-graph \hspace*{1cm}  \hfil A
 2-walk in a 4-graph}
 \caption{Three examples on an $s$-walk in a hypergraph}
 \label{2walk}
\end{figure}

Let $H$ be an $r$-uniform hypergraph with the vertex
set $V(H)$ (or $V$ for short)  and the edge set $E(H)$. We assume
$|V(H)|=n$ and $E(H) \subseteq {V\choose r}$. For
a vertex subset $S$ such that $|S|<r$, the {\it neighborhood}
$\Gamma(S)$ is $\{T | S\cap T=\emptyset \mbox { and } S\cup T
\mbox{ is an edge in } H\}$.
Let the {\it degree $d_S$} of $S$ in $H$ be the number of edges
containing $S$, i.e, $d_S=|\Gamma(S)|$.  For $1\leq s\leq r-1$, an $s$-walk of length $k$ is
a sequence of vertices
 $$v_1, v_2,\ldots, v_j,\ldots, v_{(r-s)(k-1)+r}$$
 together with
a sequence of edges $F_1, F_2, \ldots, F_k$
such that
$$F_i=\{v_{(r-s)(i-1)+1},v_{(r-s)(i-1)+2}, \ldots, v_{(r-s)(i-1)+r}\}$$
for $1\leq i\leq k$.
Here are some examples of $s$-walks as shown in Figure \ref{2walk}.

For each $i$ in $\{0,1,\ldots,k\}$, the {\it $i$-th stop} $x_i$  of the $s$-walk
is the ordered $s$-tuple \linebreak
$(v_{(r-s)i+1},v_{(r-s)i+2}, \ldots, v_{(r-s)i+s})$.
 The initial stop is $x_0$, and the terminal stop is  $x_k$.
An $s$-walk is called an $s$-path if every stop (as an ordered
$s$-tuple) is different from each other. If $x_0=x_k$, then an $s$-walk
is {\it closed}. An $s$-cycle is a closed $s$-path.

 For $1\leq s\leq r-1$ and $x,y \in \vs$, the $s$-distance
$d^{(s)}(x,y)$ is the minimum integer $k$ such that there exists an
$s$-path of length $k$ starting from $x$ and ending at $y$.  A hypergraph $H$ is
 {\it $s$-connected} if $d^{(s)}(x,y)$ is finite for every pair $(x,y)$.
If $H$ is $s$-connected, then the {\it $s$-diameter} of $H$ is
 the maximum value of $d^{(s)}(x,y)$ for $x,y\in \vs$.
%

A random $s$-walk with initial stop $x_0$   is an $s$-walk generated
as follows. Let $x_0$ be the sequence of visited vertices at initial
step. At each step, let $S$ be the set of last $s$ vertices in the
sequence of visited vertices. A random $(r-s)$-set $T$ is chosen
from  $\Gamma(S)$  uniformly; the vertex in $T$ is added into the
sequence one by one in an arbitrary order.

For $0\leq \alpha\leq 1$, an $\alpha$-lazy random $s$-walk is a
modified random $s$-walk such that with probability $\alpha$,  one can
stay at the current stop; with probability $1-\alpha$, append
$r-s$ vertices to the sequence as selected in a random $s$-walk.

For $x\in \vs$, let $[x]$ be the $s$-set consisting of
the coordinates of  $x$.

\subsection{Case $1\leq s\leq r/2$}
For $1\leq s\leq r/2$, we define a weighted undirected graph
$G^{(s)}$ over the vertex set $\vs$ as follows.
Let the weight $w(x,y)$ be $|\{F \in E(H): [x] \sqcup [y] \subseteq
F\}|$. Here $[x]\sqcup [y]$ is the disjoint union of $[x]$ and $[y]$.
In particular, if $[x]\cap [y]\not=\emptyset$, then $w(x,y)=0$.

For $x\in \vs$, the degree of $x$ in $G^{(s)}$, denoted by
$d_x^{(s)}$, is given by
\begin{equation}\label{eq:4}
d_x^{(s)}=\sum_y w(x,y)=d_{[x]}\binom{r-s}{s}s!.
 \end{equation}
Here $d_{[x]}$ means the degree of the set $[x]$ in the hypergraph
$H$. When we restrict  an $s$-walk  on $H$ to its stops, we  get a
walk  on $G^{(s)}$. This restriction keeps the length of the walk.
Therefore, the $s$-distance $d^{(s)}(x,y)$ in $H$ is simply the
graph distance between $x$ and $y$ in $G^{(s)}$; the $s$-diameter of
$H$ is simply the diameter of the graph $G^{(s)}$.

A random $s$-walk on $H$ is essentially a random walk on $G^{(s)}$.
It can be constructed from a random walk on $G^{(s)}$ by inserting
additional random $r-2s$ vertices $T_i$ between two consecutive
stops $x_i$ and $x_{i+1}$ at time $i$, where $T_i$ is chosen uniformly from
$\Gamma([x_i]\cup [x_{i+1}])$ and inserted between $x_i$ and
$x_{i+1}$ in an arbitrary order.

 Therefore, we  define the {\it $s$-th Laplacian} $\L^{(s)}$ of $H$ to be
the Laplacian of the weighted  undirected graph $G^{(s)}$.
%

%

The eigenvalues of  $\L^{(s)}$ are listed as $\lambda^{(s)}_0,
\lambda^{(s)}_1, \ldots, \lambda^{(s)}_{{n\choose s}s!}$ in the
non-decreasing order. Let
$\lambda_{\max}^{(s)}=\lambda^{(s)}_{{n\choose s}s!}$ and $\bar
\lambda^{(s)}=\max\{|1-\lambda_1^{(s)}|,|1-\lambda_{\max}^{(s)}|\}$.
For some hypergraphs, the numerical values of $\lambda_1^{(s)}$ and
$\lambda_{\max}^{(s)}$ are shown in Table \ref{tab:1} at the end of
this section.

\subsection{The case $r/2<s\leq r-1$}
For $r/2<s\leq r-1$,  we define a directed graph $D^{(s)}$ over the
vertex set $\vs$ as follows. For $x,y\in \vs$ such that
$x=(x_1,\ldots, x_s)$ and $y=(y_1,\ldots, y_s)$, let $(x,y)$ be a
directed edge  if $x_{r-s+j}=y_j$ for $1\leq j\leq 2s-r$ and
$[x]\cup [y]$ is an edge of $H$.

For $x\in\vs$, the out-degree $d^+_x$ in $D^{(s)}$  and the in-degree $d^-_x$ in $D^{(s)}$
satisfy
$$d^+_x=d_{[x]}(r-s)!=d^-_x.$$
Thus $D^{(s)}$ is a Eulerian directed graph. We write $d^{(s)}_x$ for both $d^+_x$
and $d^-_x$. Now $D^{(s)}$ is
strongly connected if and only if it is weakly connected.

Note that an $s$-walk on $H$ can be naturally viewed as a walk on
$D^{(s)}$ and vice versa. Thus the $s$-distance $d^{(s)}(x,y)$ in
$H$ is exactly the directed distance from $x$ to $y$ in $G^{(s)}$;
the $s$-diameter of $H$ is the diameter of $D^{(s)}$. A random
$s$-walk on $H$ is one-to-one corresponding to a random walk on
$D^{(s)}$.



For $\frac{r}{2}<s\leq r-1$, we define the $s$-th Laplacian
$\L^{(s)}$ as the Laplacian of the Eulerian directed graph $D^{(s)}$
(see section 2).

The eigenvalues of  $\L^{(s)}$ are listed as $\lambda^{(s)}_0,
\lambda^{(s)}_1, \ldots, \lambda^{(s)}_{{n\choose s}s!}$ in the
non-decreasing order. Let
$\lambda_{\max}^{(s)}=\lambda^{(s)}_{{n\choose s}s!}$ and $\bar
\lambda^{(s)}=\max\{|1-\lambda_1^{(s)}|,|1-\lambda_{\max}^{(s)}|\}$.
For some hypergraphs, the numerical values of $\lambda_1^{(s)}$ and
$\lambda_{\max}^{(s)}$ are shown in Table \ref{tab:1} at the end of
this section.

\subsection{Examples}
Let $K^r_n$ be the complete $r$-uniform hypergraph on $n$ vertices.
Here we compute the values of $\lambda_1^{(s)}$ and
$\lambda_{\max}^{(s)}$ for some $K^r_n$ (see Table \ref{tab:1}).

\begin{table}[htb]
\centering
\begin{tabular}{|c|c|c|c|c||c|c|c|c|}
\hline $H$  & $\lambda_1^{(4)}$ & $\lambda_1^{(3)}$  &
$\lambda_1^{(2)}$& $\lambda_1^{(1)}$
 &$\lambda_{\max}^{(1)}$   & $\lambda_{\max}^{(2)}$ & $\lambda_{\max}^{(3)}$ & $\lambda_{\max}^{(4)}$ \\
\hline $K_6^3$ & &  & $3/4$ & $6/5$ & $6/5$  & $3/2$ &   &  \\
\hline $K_7^3$   & & & $7/10$ & $7/6$
 & 7/6  & $3/2$  &  &  \\
\hline $K_6^4$ & &  $1/3$ &  $5/6$  & $6/5$
 &  $6/5$  & $3/2$ &  $1.76759$ & \\ \hline
 $K_7^4$   & &  $3/8$ &  $9/10$ & $7/6$
& $7/6$ & $7/5$ &  $7/4$ & \\
\hline $K_6^5$   & $0.1464$  & $1/2$ &  $5/6$  & $6/5$
& $6/5$ & $3/2$ & $3/2$ & $1.809$ \\
\hline $K_7^5$ & $0.1977$ & $5/8$ & $9/10$  & $7/6$ & $7/6$ & $7/5$
& $3/2$ & $1.809$ \\ \hline
\end{tabular}
\caption{The values of $\lambda_{1}^{(s)}$ and
$\lambda_{\max}^{(s)}$ of some complete hypergraphs $K_n^r$.}
\label{tab:1}
\end{table}


{\noindent}
{\bf Remark:}
From the table above, we observe $\lambda_1^{(s)}=\lambda_{\max}^{(s)}$
for some complete hypergraphs. In fact, this is true for any complete hypergraph $K_n^r$.
We point out the following fact without proofs.
For an $r$-uniform hypergraph $H$ and an integer $s$ such that
$1\leq s\leq \frac{r}{2}$, $\lambda_1^{(s)}(H)=\lambda_{\max}^{(s)}(H)$ holds if and only if $s=1$
and $H$ is a $2$-design.

\section{Properties of Laplacians}
In this section, we prove some properties of the Laplacians for
hypergraphs.
\begin{lemma}\label{lemma2}
For $1\leq s\leq r/2$, we have the following properties.
\begin{enumerate}
\item The $s$-th Laplacian has ${n\choose s}s!$ eigenvalues and all of them are in $[0,2]$.
\item The number of $0$ eigenvalues is the number of connected components
in $G^{(s)}$.
\item The Laplacian $L^{(s)}$ has an eigenvalue $2$ if and only if $r=2s$ and
$G^{(s)}$ has a bipartite component.
\end{enumerate}
\end{lemma} {\bf Proof:} Items 1 and 2 follow  from the facts
of the Laplacian of $G^{(s)}$.  If
$\L^{(s)}$ has an eigenvalue 2, then $G^{(s)}$ has a bipartite
component $T$. We want to show $r=2s$.  Suppose $r \geq 2s+1$.
Let $\{v_0,v_2,\ldots, v_{r-1}\}$ be an edge in $T$. For $0\leq i\leq 2s$ and $0\leq j\leq s-1$,
let $g(i,j)=is+j \mod (2s+1)$ and $x_i=(v_{g(i,0)},\ldots, v_{g(i,s-1)})$. Observe
$x_0, x_1,\ldots, x_{2s}$ form an odd cycle in $G^{(s)}$.
Contradiction. \hfill
$\square$

The following lemma compares $\lambda^{(s)}_1$ and
$\lambda^{(s)}_{\max}$ for different $s$.
\begin{lemma}\label{lemma3}
Suppose that $H$ is an $r$-uniform hypergraph. We have
\begin{eqnarray}
\label{labmda1} \lambda_1^{(1)} \geq \lambda_1^{(2)} \geq \ldots
\geq
\lambda_1^{(\lfloor r/2 \rfloor)};\\
\label{lambdamax} \lambda_{\max}^{(1)} \leq \lambda_{\max}^{(2)}
\leq \ldots \leq \lambda_{\max}^{(\lfloor r/2 \rfloor)}.
\end{eqnarray}
\end{lemma}
{\bf Remark:} We do not know whether similar inequalities hold for $s>\frac{r}{2}$.\\
{\bf Proof:} Let $T_s$ be the diagonal
matrix of degrees in $G^{(s)}$ and $R^{(s)}(f) $ be the Rayleigh
quotient of $\L^{(s)}$.
It suffices to show $\lambda_1^{(s)} \leq \lambda_1^{(s-1)}$ for $2 \leq  s \leq r/2$.
Recall that
$\lambda_1^{(s)}$ can be defined via the Rayleigh quotient, see
equation (\ref{eq:3}). Pick a function $f: V^{\underline{(s-1)}}
\rightarrow R$ such that $\langle f, T_{s-1} \one \rangle=0$ and
$\lambda_1^{(s-1)}=R^{(s-1)}(f)$. We define $g: \vs \rightarrow
R$ as follows
$$
g(x)=f(x'),$$
where $x'$ is a $(s-1)$-tuple consisting of the first $(s-1)$
coordinates of $x$ with the same order in $x$.
Applying equation \ref{eq:4}, we get
 $$\langle g, T_{s} \one \rangle= \sum_{x \in \vs} d_x^{(s)}g(x)=\sum_{x \in \vs}  g(x) d_{[x]}
\binom{r-s}{s}s!.$$ We have
\begin{eqnarray*}
\sum_{x }  g(x) d_{[x]} &=&\sum_{x } \sum_{F: [x] \subseteq
F} g(x)\\
&=&\sum_{x'} \sum_{F: [x'] \subseteq F} (r-s+1) f(x')\\
&=&\sum_{x'}
d_{[x']}(r-s+1) f(x')\\
&=& \frac{r-s+1}{\binom{r-s+1}{s-1}(s-1)!} \sum_{x'}
f(x')d_{x'}^{(s-1)}=0.
 \end{eqnarray*}
Here the second last equality follows from equation \ref{eq:4} and
the last one follows from the choice of $f$. Therefore,
$$\sum_x g(x)d_x^{(s)}=(r-s+2)(r-s+1)\sum_{x'}f(x')d_{x'}^{(s-1)}.$$
Thus $\langle g, T_s \one \rangle =0$. Similarly, we have
$$
\sum_x g(x)^2d_x^{(s)}=(r-s+2)(r-s+1)\sum_{x'}f(x')^2d_{x'}^{(s-1)}.
$$
Putting them together, we obtain
$$
\sum_x g(x)^2
d_x^{(s)}=(r-s+2)(r-s+1)\sum_{x'}f(x')^2d_{x'}^{(s-1)}.
$$
By the similar counting method, we have
\begin{eqnarray*}
\sum_{x \sim y} (g(x)-g(y))^2 w(x,y)&=&\sum_{x \sim y} \sum_{F: [x]
\sqcup [y] \subseteq F}(g(x)-g(y))^2\\
&=&\sum_{x' \sim y'} \sum_{F: [x'] \sqcup [y'] \subseteq
F}(r-s+1)(r-s+2) (f(x')-f(y'))^2\\
&=&(r-s+1)(r-s+2)\sum_{x' \sim y'}(f(x')-f(y'))^2w(x',y').
\end{eqnarray*}
Thus, $R^{(s)}(g)=R^{(s-1)}(f)=\lambda_1^{(s-1)}$ by the choice of
$f$.  As $\lambda_1^{(s)}$ is the infimum over all $g$, we get
$\lambda_1^{(s)} \leq \lambda_1^{(s-1)}$.

The inequality (\ref{lambdamax}) can be proved in a similarly way.
Since $\lambda^{(s)}_{\max}$ is the supremum of the Raleigh quotient,
the direction of inequalities are reversed.
\hfill $\square$

\begin{lemma}\label{lemma4}
For $r/2<s\leq r-1$, we have the following facts.
\begin{enumerate}
\item The $s$-th Laplacian has ${n\choose s}s!$ eigenvalues and all of them are in $[0,2]$.
\item The number of $0$ eigenvalues is the number of strongly connected components
in $D^{(s)}$.
\item If $2$ is an eigenvalue of $L^{(s)}$, then one of the
$s$-connected components of $H$ is bipartite.
\end{enumerate}
\end{lemma}
The proof is trivial and will be omitted.
%
%

\section{Applications}
We show some applications of Laplacians $\L^{(s)}$ of hypergraphs in
this section.
\subsection{The random $s$-walks on hypergraphs}
For $0\leq \alpha<1$ and $1\leq s\leq r/2$, after restricting an
$\alpha$-lazy random $s$-walk on a hypergraph $H$ to its stops (see
section \ref{section3}), we get an $\alpha$-lazy random walk on the
corresponding
 weighted graph $G^{(s)}$. Let $\pi(x)=d_{x}/\vol(\vs)$ for any $x\in \vs$,
where $d_x$ is the degree of $x$ in $G^{(s)}$ and $\vol(\vs)$ is the
volume of $G^{(s)}$. Applying theorem \ref{thm1}, we have the
following theorem.
\begin{theorem}
For $1\leq s\leq r/2$, suppose that $H$ is an $s$-connected
$r$-uniform hypergraph $H$ and $\lambda_1^{(s)}$ (and
$\lambda_{\max}^{(s)}$) is the first non-trivial (and the last)
eigenvalue of the $s$-th Laplacian of $H$. For $0 \leq \alpha < 1$,
the joint distribution $f_k$ at the $k$-th stop of the $\alpha$-lazy
random walk at time $k$ converges  to the stationary distribution
$\pi$ in probability. In particular, we have
$$\| (f_k-\pi)T^{-1/2} \| \leq (\bar \lambda_\alpha^{(s)})^k \|(f_0-\pi)T^{-1/2}\|,$$
where $\bar
\lambda_\alpha^{(s)}=\max\{|1-(1-\alpha)\lambda_1^{(s)}|,
|(1-\alpha)\lambda_{\max}^{(s)}-1|$, and $f_0$ is the probability
distribution at the initial stop.
\end{theorem}

For $0< \alpha<1$ and $r/2<s\leq r-1$, when restricting an
$\alpha$-lazy random $s$-walk on a hypergraph $H$ to its stops (see
section 2), we get an $\alpha$-lazy random walk on the corresponding
directed graph $D^{(s)}$.  Let $\pi(x)=d_{x}/\vol(\vs)$ for any
$x\in \vs$, where $d_x$ is the degree of $x$ in $D^{(s)}$ and
$\vol(\vs)$ is the volume of $D^{(s)}$. Applying theorem \ref{thm5},
we have the following theorem.
\begin{theorem}
For $r/2<s\leq r-1$, suppose that $H$ is an $s$-connected
 $r$-uniform hypergraph and
$\lambda_1^{(s)}$ is the first non-trivial eigenvalue of the $s$-th
Laplacian of $H$. For $0 < \alpha < 1$, the joint distribution $f_k$
at the $k$-th stop of the $\alpha$-lazy random walk at time $k$
converges  to the stationary distribution $\pi$ in probability. In
particular, we have
$$\| (f_k-\pi)T^{-1/2} \| \leq (\sigma_\alpha^{(s)})^k \|(f_0-\pi)T^{-1/2}\|,$$
where $\sigma_\alpha^{(s)}\leq \sqrt{1
-2\alpha(1-\alpha)\lambda_1^{(s)}}$, and $f_0$ is the probability
distribution at the initial stop.
\end{theorem}

{\bf Remark:} The reason why we require $0 < \alpha < 1$ in the case
$ r/2 < s \leq r-1 $ is $\sigma_0(D^{(s)})=1$ for $ r/2 < s \leq
r-1$.

\subsection{The $s$-distances and $s$-diameters in  hypergraphs}

Let $H$ be an $r$-uniform hypergraph. For $1\leq s\leq r-1$ and $x,y
\in \vs$, the $s$-distance $d^{(s)}(x,y)$ is the minimum integer $k$
such that there is an $s$-path of length $k$ starting at $x$ and
ending at $y$. For $X, Y\subseteq \vs$, let $d^{(s)}(X,Y)=\min\{d^{(s)}(x,y)\mid x\in X, y\in Y\}$.
If $H$ is $s$-connected, then the $s$-diameter $\diam^{(s)}(H)$ satisfies
$$\diam^{(s)}(H)=\max_{x,y\in \vs}\{d^{(s)}(x,y)\}.$$
For $1\leq s\leq \frac{r}{2}$, the $s$-distances in $H$ (and the
$s$-diameter of $H$) are simply the graph distances in $G^{(s)}$ (and
the diameter of $G^{(s)}$), respectively. Applying Theorem \ref{thm2} and
Corollary \ref{cor1}, we have the following theorems.

\begin{theorem}
Suppose $H$ is  an $r$-uniform hypergraph.
For integer $s$ such that $1\leq s\leq
\frac{r}{2}$, let $\lambda_1^{(s)}$ (and $\lambda_{\max}^{(s)}$) be the first non-trivial
(and the last) eigenvalue of the $s$-th Laplacian of $H$. Suppose
$\lambda_{\max}^{(s)}>\lambda_1^{(s)}>0$. For
$X,Y \subseteq \vs$, if $d^{(s)}(X,Y)\geq 2$, then we have
$$
d^{(s)}(X,Y) \leq \left \lceil \frac{\log \sqrt {\frac{\vol(\bar
X)\vol(\bar Y)}{\vol(X)\vol(Y)}}}
  {\log
  \frac{\lambda_{\max}^{(s)}+\lambda_1^{(s)}}{\lambda_{\max}^{(s)}-\lambda_1^{(s)}}}\right\rceil.
$$
Here $\vol(*)$ are volumes in $G^{(s)}$.
\end{theorem}
{\bf Remark: } We know $\lambda^{(s)}_1>0$ if and only if $H$ is
$s$-connected. The condition $\lambda^{(s)}_{\max}>\lambda^{(s)}_1$
holds unless $s=1$ and every pair of vertices is covered by edges
evenly (i.e., $H$ is a $2$-design).

\begin{theorem}
Suppose $H$ is  an $r$-uniform hypergraph.
For integer $s$ such that $1\leq s\leq
\frac{r}{2}$, let $\lambda_1^{(s)}$ (and $\lambda_{\max}^{(s)}$) be the first non-trivial
(and the last) eigenvalue of the $s$-th Laplacian of $H$.
If $\lambda_{\max}^{(s)}>\lambda_1^{(s)}>0$, then the $s$-diameter of
an $r$-uniform hypergraph $H$ satisfies
$$
\textrm{diam}^{(s)}(H) \leq   \left \lceil \frac{\log
\frac{\vol(\vs)}{\delta^{(s)}}}
  {\log
  \frac{\lambda^{(s)}_{\max}+\lambda^{(s)}_1}{\lambda^{(s)}_{\max}-\lambda^{(s)}_1}}\right\rceil.
$$
Here $\vol(\vs)=\sum_{x\in \vs}d_x=|E(H)|\frac{r!}{(r-2s)!}$ and
$\delta^{(s)}$ is the minimum degree in $G^{(s)}$.
\end{theorem}

When $ r/2 < s \leq r-1$, the $s$-distances in $H$ (and the
$s$-diameter of $H$) is the directed distance in $D^{(s)}$ (and the
diameter of $D^{(s)}$), respectively. Applying Theorem
\ref{thm8} and its remark,  we have the following theorems.
\begin{theorem}
Let $H$ be an $r$-uniform hypergraph.
For  $r/2 < s \leq r-1$ and $X,Y \subseteq
\vs$, if $H$ is $s$-connected, then we have
$$
d^{(s)}(X,Y) \leq \left \lfloor \frac{\log \frac{\vol(\bar
X)\vol(\bar Y)}{\vol(X)\vol(Y)}}
  {\log \frac{2}{2-\lambda_1^{(s)}
    }}\right\rfloor+1.
$$
Here $\lambda_1^{(s)}$ is the first non-trivial eigenvalue of the
Laplacian of  $D^{(s)}$, and $\vol(*)$ are volumes in $D^{(s)}$.
\end{theorem}

\begin{theorem}
For  $r/2 < s \leq r-1$, suppose that an $r$-uniform hypergraph $H$
is $s$-connected. Let $\lambda_1^{(s)}$ be the smallest nonzero
eigenvalue of the Laplacian of  $D^{(s)}$. The $s$-diameter of
$H$ satisfies
$$
\textrm{diam}^{(s)}(H) \leq   \left \lceil \frac{2\log
\frac{\vol(\vs)}{\delta^{(s)}}}
  {\log \frac{2}{2-\lambda_1^{(s)} }}\right\rceil.
$$
Here $\vol(\vs)=\sum_{x\in \vs}d_x=|E(H)|r!$ and $\delta^{(s)}$ is
the minimum degree in $D^{(s)}$.
\end{theorem}

\subsection{The edge expansions in hypergraphs}
In this subsection, we prove some results on the edge
expansions in hypergraphs.

Let $H$ be an $r$-uniform hypergraph. For $S\subseteq {V\choose s}$, we
recall that the volume of $S$ satisfies
$$\vol(S)=\sum_{x\in S}d_x.$$
Here $d_x$ is the degree of the set $x$ in $H$.  In particular,
we have
$$\vol\left({V\choose s}\right)= |E(H)|{r\choose s}.$$
The {\it density} $e(S)$ of $S$ is $\frac{\vol(S)}{\vol({V\choose s})}$.
Let $\bar S$ be the complement set of $S$ in ${V\choose s}$. We have
$$e(\bar S)=1-e(S).$$

For $1\leq t\leq s\leq r-t$, $S\subseteq {V\choose s}$, and $T\subseteq
{V\choose t}$, let
$$E(S,T)=\{F\in E(H): \exists x\in S, \exists y\in T, x\cap y=\emptyset,
\mbox{ and } x\cup y\subseteq F\}.$$ Note that $|E(S,T)|$ counts the
number of edges contains $x\sqcup y$ for some $x\in S$ and  $y\in T$.

Particularly, we have
$$\left | E\left({V\choose s},{V\choose t}\right) \right |=|E(H)|\frac{r!}{s!t!(r-s-t)!}.$$
\begin{theorem} \label{thm11}
For $1\leq t\leq s \leq \frac{r}{2}$, $S\subseteq {V\choose s}$, and
$T\subseteq {V\choose t}$, let $e(S,T)=\frac{|E(S,T)|}{|E({V\choose
s},{V\choose t})|}$. We have
\begin{equation}
  \label{eq:exp1}
  |e(S,T)-e(S)e(T)|\leq \bar\lambda^{(s)}\sqrt{e(S)e(T)e(\bar S)e(\bar T)}.
\end{equation}
\end{theorem}
{\bf Proof:} Let $G^{(s)}$ be the weighed undirected graph defined
in section \ref{section3}. Define $S'$ and $T'$ (sets of ordered
$s$-tuples) as follows
$$S'=\{x\in \vs \mid [x]\in S\};$$
$$T'=\{(y,z)\in \vs \mid [y] \in T \}.$$
Let $\bar S'$ (or $\bar T'$) be the complement set of $S'$ (or $T'$) in $\vs$,
respectively.  We make a convention that
$\vol_{G^{(s)}}(\ast)$ denotes volumes in $G^{(s)}$ while
$\vol(\ast)$ denotes volumes $H$. We have
\begin{eqnarray}
\label{eq:start}
  \vol_{G^{(s)}}(G^{(s)})&=&\vol\left({V\choose s} \right)\frac{s!(r-s)!}{(r-2s)!}; \\
  \vol_{G^{(s)}}(S')&=&\vol(S)\frac{s!(r-s)!}{(r-2s)!};\\
  \vol_{G^{(s)}}(T')&=&\vol(T)\frac{t!(r-t)!}{(r-2s)!};\\
\vol_{G^{(s)}}(\bar S')&=&\vol(\bar S)\frac{s!(r-s)!}{(r-2s)!};\\
  \vol_{G^{(s)}}(\bar T')&=&\vol(\bar T)\frac{t!(r-t)!}{(r-2s)!}.
\label{eq:end}
\end{eqnarray}
Let $E_{G^{(s)}}(S',T')$ be the number of edges between $S'$ and
$T'$ in $G^{(s)}$. We get
$$|E_{G^{(s)}}(S',T')|=\frac{(r-s-t)! s! t!}{(r-2s)!}|E(S,T)|.$$
Applying Theorem \ref{thm3} to the sets $S'$ and $T'$ in $G^{(s)}$,
we obtain
\begin{eqnarray*}
   &&\left| |E_{G^{(s)}}(S',T')| -\frac{\vol_{G^{(s)}}(S')\vol_{G^{(s)}}(T') }{\vol_{G^{(s)}}(G^{(s)})}
\right| \hspace*{2in}\\
&& \hspace*{1in}\leq \bar\lambda_1^{(s)}\frac{
\sqrt{\vol_{G^{(s)}}(S')\vol_{G^{(s)}}(T') \vol_{G^{(s)}}(\bar
S')\vol_{G^{(s)}}(\bar T')} }{\vol_{G^{(s)}}(G^{(s)})}.
\end{eqnarray*}
Combining equations (\ref{eq:start}-\ref{eq:end}) and the inequality
above, we obtain inequality \ref{eq:exp1}. \hfill $\square$

Now we consider the case that $s>\frac{r}{2}$. Due to the fact that
$\sigma_0^{(s)}=1$, we have to use the weaker expansion theorem
\ref{thm7}. Note that
$$\left | E\left({V\choose s},{V\choose t}\right) \right |=|E(H)|\frac{r!}{(r-s-t)!s!t!}.$$
We get the following  theorem.
\begin{theorem}\label{t:exp2}
For $1\leq t<\frac{r}{2}<s< s+t\leq r$, $S\subseteq {V\choose s}$, and
$T\subseteq {V\choose t}$, let $e(S,T)=\frac{|E(S,T)|}{|E({V\choose
s},{V\choose t})|}$. If $|x\cap y| \not=\min\{t, 2s-r\}$
for any $x\in S$ and $y\in T$, then we have
\begin{equation}
  \label{eq:exp2}
  |\frac{1}{2}e(S,T)-e(S)e(T)|\leq \bar\lambda^{(s)}\sqrt{e(S)e(T)e(\bar S)e(\bar T)}.
\end{equation}
\end{theorem}




{\bf Proof:} Recall that $D^{(s)}$ is the directed graph defined in section
\ref{section3}. Let
$$S'=\{x\in \vs \mid [x]\in S\};$$
$$T'=\{(y,z)\in \vs \mid [z]\in T\}.$$
We also denote $\bar S'$ (or $\bar T'$) be the complement set of
$S'$ (or $T'$) in $\vs$,   respectively.  We use  the convention
that $\vol_{D^{(s)}}(\ast)$ denotes the  volumes in $D^{(s)}$
while $\vol(\ast)$ denotes the  volumes in the hypergraph $H$. We
have
\begin{eqnarray}
\label{eq:start2}
  \vol_{D^{(s)}}(D^{(s)})&=&\vol\left({V\choose s} \right)s!(r-s)!; \\
  \vol_{D^{(s)}}(S')&=&\vol(S)s!(r-s)!;\\
  \vol_{D^{(s)}}(T')&=&\vol(T)t!(r-t)!;\\
\vol_{D^{(s)}}(\bar S')&=&\vol(\bar S)s!(r-s)!;\\
  \vol_{D^{(s)}}(\bar T')&=&\vol(\bar T)s!(r-s)!.
\label{eq:end2}
\end{eqnarray}
Let $E_{D^{(s)}}(S',T')$ (or $E_{D^{(s)}}(T',S')$) be the number of
directed edges from $S'$ to $T'$ ( or from $T'$ to $S'$) in
$D^{(s)}$, respectively. We get
$$|E_{D^{(s)}}(S',T')|=  (r-s-t)! s! t!|E(S,T)|.$$
From the condition  $|x\cap y|\not=\min\{t,2s-r\}$ for each $x\in S$
and each $y\in T$, we observe
$$ E_{D^{(s)}}(T',S')=0.$$

Applying Theorem \ref{thm7} to the sets $S'$ and $T'$ in $D^{(s)}$,
we obtain
\begin{eqnarray*}
   &&\left| \frac{|E_{D^{(s)}}(S',T')|+|E_{D^{(s)}}(T',S')|}{2} -\frac{\vol_{D^{(s)}}(S')\vol_{D^{(s)}}(T') }{\vol_{D^{(s)}}(D^{(s)})}
\right| \hspace*{2in}\\
&& \hspace*{1in}\leq \bar\lambda_1^{(s)}\frac{
\sqrt{\vol_{D^{(s)}}(S')\vol_{D^{(s)}}(T') \vol_{D^{(s)}}(\bar
S')\vol_{D^{(s)}}(\bar T')} }{\vol_{D^{(s)}}(D^{(s)})}.
\end{eqnarray*}
Combining equations (\ref{eq:start2}-\ref{eq:end2}) and the
inequality above, we get inequality \ref{eq:exp2}.
\hfill $\square$ \\

Nevertheless, we have the following strong edge expansion theorem
for $\frac{r}{2}< s\leq r-1$. For  $S, T \subseteq {V\choose s}$, let
$E'(S,T)$ be the set of edges of the form $x\cup y$ for some $x\in
S$ and $y\in T$. Namely,
$$E'(S,T)=\{F\in E(H)\mid \exists x\in S, \exists y\in T,
F=x\cup y\}.$$  Observe that
$$\left |E'\left({V\choose s},{V\choose s}\right)\right|=|E(H)|\frac{r!}{(r-s)!(2s-r)!(r-s)!}.$$
\begin{theorem}\label{t:exp3}
For $\frac{r}{2}<s\leq r-1$ and $S, T \subseteq {V\choose s}$, let
$e'(S,T)=\frac{|E'(S,T)|}{|E'({V\choose s},{V\choose s})|}$. We have
\begin{equation}
  \label{eq:exp3}
  |e'(S,T)-e(S)e(T)|\leq \bar\lambda^{(s)}\sqrt{e(S)e(T)e(\bar S)e(\bar T)}.
\end{equation}
\end{theorem}
{\bf Proof:}
 Let
$$S'=\{x\in \vs \mid [x]\in S\};$$
$$T'=\{y\in \vs \mid [y[\in T\}.$$
Let $\bar S'$ (or $\bar T'$) be the complement set of $S'$ (or $T'$
respectively) in $\vs$.  We use  the convention that
$\vol_{D^{(s)}}(\ast)$ denotes the  volumes in $D^{(s)}$ while
$\vol(\ast)$ denotes the  volumes in the hypergraph $H$. We have
\begin{eqnarray}
\label{eq:start3}
  \vol_{D^{(s)}}(D^{(s)})&=&\vol\left({V\choose s} \right)s!(r-s)!; \\
  \vol_{D^{(s)}}(S')&=&\vol(S)s!(r-s)!;\\
  \vol_{D^{(s)}}(T')&=&\vol(T)s!(r-s)!;\\
\vol_{D^{(s)}}(\bar S')&=&\vol(\bar S)s!(r-s)!;\\
  \vol_{D^{(s)}}(\bar T')&=&\vol(\bar T)s!(r-s)!.
\label{eq:end3}
\end{eqnarray}
Let $E_{D^{(s)}}(S',T')$ (or $E_{D^{(s)}}(T',S')$) be the number of
directed edges from $S'$ to $T'$ ( or from $T'$ to $S'$)  in $D^{(s)}$,
respectively. We get
$$|E_{D^{(s)}}(S',T')|= |E_{D^{(s)}}(T',S')|=(r-s)! (2s-r)! (r-s)!|E'(S,T)|.$$
Applying Theorem \ref{thm7} to the sets $S'$ and $T'$ on $D^{(s)}$,
we obtain
\begin{eqnarray*}
   &&\left| \frac{|E_{D^{(s)}}(S',T')|+|E_{D^{(s)}}(T',S')|}{2} -\frac{\vol_{D^{(s)}}(S')\vol_{D^{(s)}}(T') }{\vol_{D^{(s)}}(D^{(s)})}
\right| \hspace*{2in}\\
&& \hspace*{1in}\leq \bar\lambda_1^{(s)}\frac{
\sqrt{\vol_{D^{(s)}}(S')\vol_{D^{(s)}}(T') \vol_{D^{(s)}}(\bar
S')\vol_{D^{(s)}}(\bar T')} }{\vol_{D^{(s)}}(D^{(s)})}.
\end{eqnarray*}
Combining equations (\ref{eq:start3}-\ref{eq:end3}) and the
inequality above, we get inequality \ref{eq:exp3}. \hfill $\square$

\section{Concluding Remarks}
In this paper, we introduced a set of Laplacians for $r$-uniform
hypergraphs. For $1\leq s\leq r-1$,  the $s$-Laplacian $\L^{(s)}$ is
derived from the random $s$-walks on hypergraphs. For $1\leq s\leq
\frac{r}{2}$, the $s$-th Laplacian $\L^{(s)}$ is defined to be the
Laplacian of the corresponding weighted graph $G^{(s)}$. The first Laplacian
$\L^{(1)}$ is exactly the Laplacian introduced by  Rodr{\`i}guez
\cite{rod}.

For $\frac{r}{2}\leq s\leq r-1$, the $\L^{(s)}$ is defined to be the
Laplacian of the corresponding Eulerian directed graph $D^{(s)}$. At first glimpse,
$\sigma_0(D^{(s)})$ might be a good parameter. However, it is not
hard to show that $\sigma_0(D^{(s)})=1$ always holds, which makes
Theorem \ref{thm6} useless for hypergraphs. We can use
weaker Theorem \ref{thm7} for hypergraphs. Our work is based on
(with some improvements)
 Chung's recent work \cite{fan5,
fan6}  on directed graphs.

Let us recall Chung's definition of Laplacians \cite{fan2} for
regular hypergraphs. An $r$-uniform hypergraph $H$ is $d$-regular if
$d_x=d$ for every $x\in \vrm$. Let $G$ be a graph on the vertex set
$\vrm$. For $x,y\in \vrm$, let $xy$ be an edge if
$x=x_1x_2,\ldots,x_{r-1}$ and $y=y_1x_2,\ldots,x_{r-1}$ such that
$\{x_1,y_1, x_2,\ldots,x_{r-1}\}$ is an edge of $H$. Let $A$ be the
adjacency matrix of $G$, $T$ be the diagonal matrix of degrees in
$G$, and $K$ be the adjacency matrix of the complete graph on the
edge set $\vrm$. Chung \cite{fan2} defined the Laplacian $\L$ such
that
$$\L=T-A + \frac{d}{n}(K+(r-1)I).$$
This definition comes from the homology theory of hypergraphs.
Firstly,  $\L$ is not normalized in Chung's definition, i.e., the
eigenvalues are not in the interval $[0,2]$. Secondly, the add-on
term $\frac{d}{n}(K+(r-1)I)$ is not related to the structures of $H$.
If we ignore the add-on term and normalize the matrix, we
essentially get the Laplacian of the graph $G$. Note $G$ is
disconnected, then $\lambda_1(G)=0$ and it is not interesting. Thus
 Chung added the additional term. The graph $G$ is actually very
closed to our Eulerian directed graph $D^{(r-1)}$. Let $B$ be the
adjacency matrix of $D^{(r-1)}$. In fact we have $B=QA$, where $Q$
is a rotation  which maps $x=x_1,x_2\ldots,x_{r-1}$ to
$x'=x_2\ldots,x_{r-1},x_1$. Since $d_{x}=d_{x'}$,  $Q$ and $T$
commute, we have
\begin{eqnarray*}
  (T^{-1/2}BT^{-1/2})'(T^{-1/2}BT^{-1/2}) &=& T^{-1/2}B'T^{-1}BT^{-1/2}\\
&=&  T^{-1/2}A'Q'T^{-1}QAT^{-1/2}\\
&=& T^{-1/2}A'T^{-1}Q'QAT^{-1/2}\\
&=& T^{-1/2}A'T^{-1} AT^{-1/2}.
\end{eqnarray*}
Here we use the fact $Q'Q=I$. This identity means that the singular
values of $I-\L^{(r-1)}$ is precisely equal to $1$ minus the
Laplacian eigenvalues of the graph $G$.

Our definitions of Laplacians $\L^{(s)}$ are clearly related to the
quasi-randomness of hypergraphs.  We are very interested in this
direction. Many concepts  such as the $s$-walk, the $s$-path, the
$s$-distance, and the $s$-diameter, have their independent interest.

\end{document}